\documentclass{amsart}
\usepackage{amssymb,amsmath}
\newtheorem{theorem}{Theorem}
\newtheorem{lemma}{Lemma}
\newtheorem{problem}{Problem}
\newcommand{\bt}{\begin{theorem}}
\newcommand{\et}{\end{theorem}}
\newcommand{\bl}{\begin{lemma}}
\newcommand{\el}{\end{lemma}}
\newcommand{\bp}{\begin{problem}}
\newcommand{\ep}{\end{problem}}
\newcommand{\pf}{{\bf Proof}.\ }
\newcommand{\eop}{$\square$\vspace{.8cm}}
\newcommand{\bal}{\begin{align*}}
\newcommand{\eal}{\end{align*}}
\newcommand{\bq}{\begin{eqnarray*}}
\newcommand{\eq}{\end{eqnarray*}}
\newcommand{\be}{\begin{eqnarray}}
\newcommand{\ee}{\end{eqnarray}}
\newcommand{\beq}{\begin{equation}}
\newcommand{\eeq}{\end{equation}}
\newcommand{\benum}{\begin{enumerate}}
\newcommand{\eenum}{\end{enumerate}}
\newcommand{\ba}{\begin{array}}
\newcommand{\ea}{\end{array}}

\begin{document}
\title{Arithmetic progressions in sets with small sumsets}
\author{J\'ozsef Solymosi}
\thanks{This research was supported by NSERC and OTKA grants.}
\address{Department of Mathematics\\
University of British Columbia\\
Vancouver} \email{solymosi@math.ubc.ca}

\begin{abstract}
We present an elementary proof that if $A$ is a finite set of
numbers, and the sumset $A+_GA$ is small, $|A+_GA|\leq c|A|$,
along a dense graph $G$, then $A$ contains $k$-term arithmetic
progressions.
\end{abstract}

\maketitle

\section{Introduction}
A well known theorem of Szemer\'edi \cite{SZ} states that every
dense subset of integers contains long arithmetic progressions. A
different, but somehow related result of Freiman \cite{FR} says
that if the sumset of a finite set of numbers $A$ is small, i.e.
$|A+A|\leq C|A|,$ then $A$ is the subset of a (not very large)
generalized arithmetic progression. Balog and Szemer\'edi proved
in \cite{BSZ} that a similar structural statement holds under
weaker assumptions. (For correct statements and details, see
\cite{NA}). As a corollary of their result, Freiman's
theorem, and Szemer\'edi's theorem about $k$-term arithmetic
progressions, Balog and Szemer\'edi proved Theorem 1 below. 
The goal of this paper is to present a simple,
purely combinatorial proof of this assertion.

Let $A$ be a set of numbers and $G$ be a graph such that the
vertex set of $G$ is $A.$ The {\em sumset of $A$ along $G$} is
\[
A+_GA = \{a+b: a,b \in A \text{ and } (a,b) \in E(G)\}.
\]

\bt For every $c,K,k >0$ there is a threshold
$n_0=n_0(c,K,k)$ such that if $|A|=n\geq n_0$, $|A+_GA|\leq K|A|$,
and $|E(G)|\geq cn^2$, then $A$ contains a $k$-term arithmetic
progression. \et

\section{Lines and hyperplanes}
There are arrangements of $n$ lines on the Euclidean plane such that 
the maximum number of points incident with at least three lines is 
${n^2\over 6}.$ Not much is known about the
structure of arrangements where the number of such points is close
to the maximum, say $cn^2$, where $c$ is a positive constant. 
Nevertheless, the following is true.

\bl        \label{lemma:line} For every $c>0$ there is a threshold
$n_0=n_0(c)$ and a positive $\delta =\delta (c)$ such that, for any
set of $n\geq n_0$ lines $L$ and any set of $m\geq cn^2$ points
$P$, if every point is incident to three lines, then there are at
least $\delta n^3$ triangles in the arrangement. (A triangle is a
set of three distinct points from $P$ such that any two are
incident to a line from $L.$) \el

\pf This lemma follows from the following theorem of Ruzsa and
Szemer\'edi \cite{RSZ}.

\bt \cite{RSZ} Let $G$ be a graph on $n$ vertices. If $G$ is the
union of $cn^2$ edge-disjoint triangles, then $G$ contains at
least $\delta n^3$ triangles, where $\delta$ depends on $c$ only.
\et

To prove Lemma 1, let us construct a graph where $L$ is the vertex
set, and two vertices are adjacent if and only if the
corresponding lines cross at a point of $P$. This graph is the
union of $cn^2$ disjoint triangles, every point of $P$ defines a
unique triangle, so we can apply Theorem 2.\eop

The result above suffices to prove Theorem 1 for 3-term
arithmetic progressions. But for larger values of $k$, we need a
generalization of Lemma 1.

\bl        \label{lemma:plane} For every $c>0$ and $d\geq 2$, there
is a threshold $n_0=n_0(c,d)$ and a positive $\delta =\delta
(c,d)$ such that, for any set of $n\geq n_0$ hyperplanes $L$ and
any set of $m\geq cn^d$ points $P$, if every point is incident to
$d+1$ hyperplanes, then there are at least $\delta n^{d+1}$
simplices in the arrangement. (A simplex is a set of $d+1$
distinct points from $P$ such that any $d$ are incident to a
hyperplane from $L.$) \el

Lemma 2 follows from the Frankl-R\"odl conjecture \cite{FRR}, the 
generalization of Theorem 2. The
$d=3$ case was proved in \cite{FRR} and the conjecture has been
proved recently by Gowers \cite{GO} and independently by Nagle,
R\"odl, Schacht, and Skokan \cite{NRS},\cite{RS}. For details, how
Lemma 2 follows from the Frankl-R\"odl conjecture, see \cite{SO}.

\section{The $k=3$ case}
Let $A$ be a set of numbers and $G$ be a graph such that the vertex
set of $G$ is $A.$ We define the {\em difference-set of $A$ along
$G$} as

\[
A-_GA = \{a-b: a,b \in A \text{ and } (a,b) \in E(G)\}.
\]

\bl For every $\epsilon ,c,K >0$ there is a number $D=D(\epsilon
,c,K)$ such that if $|A+_GA|\leq K|A|$ and $|E(G)|\geq c|A|^2$,
then there is a graph $G'\subset G$ such that $|E(G')|\geq
(1-\epsilon)|E(G)|$ and $|A-_{G'}A|\leq D|A|$. \el

\pf Let us consider the arrangement of points given by a subset of
the Cartesian product $A\times A$ and the lines $y=a$, $x=a$ for
every $a\in A$, and $x+y=t$ for every $t\in A+_GA.$ The pointset
$P$ is defined by $(a,b)\in P$ iff $(a,b)\in E(G).$ By Lemma 1,
the number of triangles in this arrangement is $\delta n^3.$ The
triangles here are right isosceles triangles. We say that a
point in $P$ is {\em popular} if the point is the right-angle
vertex of at least $\alpha n$ triangles. Selecting
$\alpha={{\delta (\epsilon c)}\over {\epsilon c}}$, where $\delta
(\cdot)$ is the function from Lemma 1, all but at most $\epsilon
cn^2$ points of $P$ are popular.

A $t\in A-A$ is {\em popular} if $|\{(a,b):a-b=t; a,b\in A\}|\geq
\alpha n.$ The number of popular $t$s is at most $Dn$, where $D$
depends on $\alpha$ only. $A\times A$ is a Cartesian product,
therefore every triangle can be extended to a square adding one
extra point from $A\times A$. Every popular point $p$ is the
right-angle vertex of at least $\alpha n$ triangles. Therefore $p$
is incident to a line $x-y=t$, where $t$ is popular, because this
line contains at least $\alpha n$ ``fourth" vertices of squares
with $p$. \eop

{\bf Proof of Theorem 1, case $k=3.$} Let us apply Lemma 1 to the
pointset $P'$ defined by $(a,b)\in P'$ iff $(a,b)\in E(G')$ and
the lines are $y=a$ for every $a\in A$, $x-y=t$ for every $t\in
A-_{G'}A$, and $x+y=s$ for every $s\in A+_GA.$ By Theorem 2, if
$|A|$ is large enough, then there are triangles in the
arrangement. The vertices of such triangles are vertices from
$P'\subset A\times A.$ The vertical lines through the vertices
form a 3-term arithmetic progression and therefore $A$ contains 
$\delta n^2$ 3-term arithmetic progressions, where $\delta > 0$ 
depends on $c$ only. \eop

\section{The general, $k>3$, case}

Following the steps of the proof for $k=3$, we prove the general
case by induction on $k.$ We prove the following theorem, which
was conjectured by Erd\H os and proved by Balog and Szemer\'edi in
\cite{BSZ}. Theorem 3, together with the $k=3$ case, gives a proof
of Theorem 1.

\bt For every $c>0$ and $k>3$ there is an $n_0$ such that, if $A$
contains at least $c|A|^2$ 3-term arithmetic progressions and
$|A|\geq n_0$, then $A$ contains a $k$-term arithmetic progression.
\et

Instead of triangles, we must consider simplices. Set $k=d$. 
In the $d$-dimensional space we show that $A\times \cdots \times A$,
the $d$-fold Cartesian product of $A$, contains a simplex in which the
vertices' first coordinates form a $(d+1)$-term arithmetic
progression.

The simplices we are looking for are homothetic\footnote{Here we
say that two simplices are homothetic if the corresponding facets
are parallel.} images of the simplex $S_d$ whose vertices are listed 
below:
\begin{displaymath}
\begin{array}{c}
(0, 0, 0,0, \ldots ,0,0)\\
(1, 1,0,0, \ldots ,0,0)\\
(2, 0, 1,0, \ldots ,0,0)\\
(3, 0,0,1, \ldots ,0,0)\\
\vdots\\
(d-1, 0, \ldots ,1,0)\\
(d, 0,0,0, \ldots ,0,0).
\end{array}
\end{displaymath}
\indent An important property of $S_d$ is that its facets can be pushed into
a ``shorter" grid. The facets of $S_d$ are parallel to hyperplanes,
defined by the origin $(0,0,0,0,\ldots ,0,0)$, and some
$(d-1)$-tuples of the grid
$$\{0,1,2,\ldots ,d-1\}\times \{-1,0,1\}\times\{0,1\}^{d-2}.$$

For example, if $d=3$, then the facets are
\begin{displaymath}
\begin{array}{c}
\{(0,0,0),(1,1,0),(2,0,1)\}\\
\{(0,0,0),(1,1,0),(3,0,0)\}\\
\{(0,0,0),(2,0,1),(3,0,0)\}\\
\{(1,1,0),(2,0,1),(3,0,0)\},
\end{array}
\end{displaymath}
and the corresponding parallel planes in
$$\{0,1,2\}\times
\{-1,0,1\}\times\{0,1\}$$ are the planes incident to the triples
\begin{displaymath}
\begin{array}{c}
\{(0,0,0),(1,1,0),(2,0,1)\}\\
\{(0,0,0),(1,1,0),(2,0,0)\}\\
\{(0,0,0),(2,0,1),(2,0,0)\}\\
\{(0,0,0),(1,-1,1),(2,-1,0)\}.
\end{array}
\end{displaymath}
\indent In general, if a facet of $S_d$ contains the origin and the ``last
point" $(d, 0,0,0, \ldots ,0,0),$ then if we replace the later one by 
$(d-1, 0,0,0, \ldots ,0,0)$, the new $d$-tuples define the same hyperplane. 
The remaining facet $f$, given by
\begin{displaymath}
\begin{array}{c}
(1, 1,0,0, \ldots ,0,0)\\
(2, 0, 1,0, \ldots ,0,0)\\
(3, 0,0,1, \ldots ,0,0)\\
\vdots\\
(d-1, 0, \ldots ,1,0)\\
(d, 0,0,0, \ldots ,0,0),
\end{array}
\end{displaymath}
is parallel to the hyperplane through the vertices of $f-(1,1,0,0,
\ldots ,0,0),$
\begin{displaymath}
\begin{array}{c}
(0, 0,0,0, \ldots ,0,0)\\
(1, -1, 1,0, \ldots ,0,0)\\
(2, -1,0,1, \ldots ,0,0)\\
\vdots\\
(d-2, -1, \ldots ,1,0)\\
(d-1,-1,0,0, \ldots ,0,0).
\end{array}
\end{displaymath}
\indent In a homothetic copy of the grid $$T_d=\{0,1,2,\ldots ,d-1\}\times
\{-1,0,1\}\times\{0,1\}^{d-2},$$ the image of the origin is called
the {\em holder} of the grid.

As the induction hypothesis, let us suppose that Theorem 3 is true
for a $k\geq 3$ in a stronger form, providing that the number of
$k$-term arithmetic progressions in $A$ is at least $c|A|^2.$

Then the number of distinct homothetic copies of $T_d$ in
$\mathbb{A}_d=\underbrace{A\times \ldots  \times A}_d$ is at least
$c'|A|^{d+1}$ ($c'$ depends on $c$ only). Let us say that a point
$p\in \mathbb{A}_d$ is \emph{popular} if $p$ is the holder of at least
$\alpha |A|$ grids. If $p$ is popular, then for any facet of
$S_d$, $f$, the point $p$ is the element of at least $\alpha |A|$
$d$-tuples, similar and parallel to $f.$ If $\alpha$ is small
enough, then at least $\gamma |A|^d$ points of $\mathbb{A}_d$ are
popular, where $\gamma$ depends on $c$ and $\alpha$ only.

A hyperplane $H$ is \emph{$\beta$-rich} if it is incident to many points,
$|H\cap \mathbb{A}_d|\geq \beta |A|^{d-1}.$ For every facet of
$S_d$, $f$, let us denote the set of $\beta$-rich hyperplanes
which are parallel to $f$ by $\mathcal{H}_f.$

\bl For some choice of $\beta$, at least half of the
popular points are incident to $d+1$ $\beta$-rich hyperplanes,
parallel to the facets of $S_d.$ \el

Suppose to the contrary that for a facet $f$, more than
${\gamma\over 2d} |A|^d$ popular points are not incident to
hyperplanes of $\mathcal{H}_f.$ Then more than

\begin{equation}
\alpha |A|{\gamma\over 2d} |A|^d={\gamma\alpha\over 2d} |A|^{d+1}
\end{equation}

\noindent $d$-tuples, similar and parallel to $f$, are not covered by
$\mathcal{H}_f.$ Let us denote the hyperplanes incident to the
``uncovered" $d$-tuples by $L_1,L_2,\ldots ,L_m$, and the number of
points on the hyperplanes by $\mathcal{L}_1,\mathcal{L}_2,\ldots
,\mathcal{L}_m.$ A simple result of Elekes and Erd\H os
\cite{EE},\cite{EL} implies that hyperplanes with few points cannot 
cover many $d$-tuples.

\bt \cite{EE} The number of homothetic copies of $f$ in $L_i$ is
at most $c_d\mathcal{L}_i^{1+1/(d-1)}$, where $c_d$ depends on $d$
only. \et

The inequalities
$$\sum_{i=1}^m\mathcal{L}_i\leq |A|^d, \text{ and } \mathcal{L}_i\leq \beta |A|^{d-1}.$$
lead us to the proof of Lemma 4.

The number of $d$-tuples covered by $L_i$s is at most

$$c_d\sum_{i=1}^m\mathcal{L}_i^{1+1/(d-1)}\leq c_d{{|A|^d}\over 
{\beta |A|^{d-1}}}(\beta |A|^{d-1})^{1+1/(d-1)}=c_d\beta^{1/(d-1)}|A|^{d+1}.$$

If we compare this bound to (1), and choose $\beta$ 
such that ${\gamma\alpha\over 2d}=c_d\beta^{1/(d-1)}$,
then at least half of the popular points are covered by $d+1$
$\beta$-rich hyperplanes parallel to the facets of $S_d.$
\eop

Finally we can apply Lemma 2 with the pointset $P$ of ``well-covered"
popular points of $\mathbb{A}_d$ and with the sets of hyperplanes
$L=\bigcup_{f\subset S_d}\mathcal{H}_f.$ The number of points is
at least ${\gamma\alpha\over 2} |A|^{d}$. For a given $f,$
$|\mathcal{H}_f|\leq {{|A|^d}\over{\beta |A|^{d-1}}}=|A|/\beta.$
The number of hyperplanes in $L$ is at most $(d+1)|A|/\beta.$ By
Lemma 2, we have at least $\delta '|A|^{d+1}$ homothetic copies of
$S_d$ in $\mathbb{A}_d.$ Let us project them onto $x_1$, the first
coordinate axis. Every image is a $(k+1)$-term arithmetic
progression, and the multiplicity of one image is at most
$|A|^{d-1}.$ Therefore there are at least $\delta '|A|^2$
$(k+1)$-term arithmetic progressions in $A.$
\eop

\section{$G_n=K_n$}
When the full sumset $A+A$ is small then it is easier to prove
that $A$ contains long arithmetic progressions. We can use the
following Pl\"unecke type inequality \cite{PL,RU,NA}.

\bt

Let $A$ and $B$ be finite subsets of an abelian group such that
$|A|=n$ and $|A+B|\leq \delta n$. Let $k\geq 1$ and $l\geq 1.$
Then

$$|kB-lB|\leq \delta^{k+l}n.$$

\et

It follows from the inequality, that for any dimension $d$ and
$d$-dimensional integer vector $\vec{v}=(x_1,\ldots ,x_d), x_i\in
\mathbb{Z}$, there is a $c>0$ depending on $d,\delta$ and $\vec{v}$
such that the following holds: \emph{If $|A+A|\leq \delta |A|$, then
$\mathbb{A}_d$ can be covered by $c|A|$ hyperplanes with the same
normal vector $\vec{v}$}. Using this, we can define our
hyperplane-point arrangement, with the hyperplanes parallel to the
facets of $S_d$ containing at least one point of $\mathbb{A}_d$,
and the pointset of the arrangement is $\mathbb{A}_d.$ Then we
do not have to deal with rich planes and popular points, and we can
apply Lemma 2 directly.

\end{document}